\documentclass[]{article}
\usepackage{amssymb}

\def\C{\mathbb{C}} 
\def\G{\mathbb{G}} 

\def\R{\mathbb{R}}  

\def\no{\noindent}

\def\beq{\begin{equation}}
\def\eeq{\end{equation}}

\def\w{\wedge}

\usepackage{graphicx}   
\usepackage{makeidx}
\def\no{\noindent}
\newtheorem{lem1}{Lemma}
\newtheorem{def1}{Definition}
\newtheorem{lem2}[lem1]{Lemma}
\newtheorem{thm1}{Theorem (Pl\"ucker's Relations)}
\newtheorem{thm2}[thm1]{Theorem}

\title{Notes on Pl\"ucker's Relations in Geometric  Algebra}
\author{Garret Sobczyk \\
	Universidad de las Am\'ericas-Puebla \\
	Departamento de F\'isico-Matem\'aticas \\
	72820 Puebla, Pue., M\'exico}
\begin{document}
	
	\maketitle
	
	\begin{abstract}Grassmannians are of fundamental importance in projective geometry, algebraic geometry, and representation theory. A vast literature has grown up utilizing using many different languages of higher mathematics, such as multilinear and tensor algebra, matroid theory, and Lie groups and Lie algebras. Here we explore the basic idea of the Pl\"ucker relations in Clifford's geometric algebra. We discover that the Pl\"ucker Relations can be fully characterized in terms of the geometric product, without the need for 
	 a confusing hodgepodge of many different formalisms and mathematical traditions found in the literature.   \end{abstract}
\smallskip
\no {\em AMS Subject Classification:} 14M15, 14A25, 15A66, 15A75

\smallskip
\no {\em Keywords: Clifford geometric algebras, Grassmannians, matroids, Pl\"ucker relations}
\section{Pl\"ucker's relations in geometric algebra}
 Clifford's {\it geometric algebra} is considered here to be the natural {\it geometrization} of the real and complex number sytems; an overview is provided in \cite{geoS2017}. We use exclusively the language of geometric algebra, developed in the books \cite{H/S,SNF}. Pl\"ucker's coordinates are briefly touched upon on page 30 of \cite{H/S}, but Pl\"ucker's relations are not discussed. These relations characterize when an $r$-vector $B$ in the geometric algebra $\G_n$ of an $n$-dimensional (real or complex) Euclidean vector space $\R^n$ (or $\C^n$), can be expressed as the outer product of $r$-vectors. Several references to the relevant literature, and its more advanced applications, are \cite{ardila,dandw91}, 
 \cite[pp.227-231]{FH1991}, \cite{harris1992}, \cite[pp.144-148]{milne}, and \cite{ctc-wall2008}. 

 Let $B=\sum_J B_{J} e_J$ be an $r$-vector expanded in the standard orthonormal basis of $r$-vector blades $e_{J_m} \in \G_n^r$, ordered lexicographically, 
 and let $\# (J_p\cap J_m)$ be the number of integers that the sequences $J_p$ and $J_m$ have in common. Given an $r$-vector $B \in \G_{n}^r$, the {\it rank space} $V_B$ of $B$ is 
 \[  V_B := \{ x \in V^n | \ x \w B = 0  \}\subset V^n,   \]
 where $V^n := V_{e_{1 \cdots n}}\equiv \R^n $ is the underlying vector space of
 the geometric algebra $\G_n$. The {\it dimension} of the subspace $V_B$ of $V^n$  
 defined by $B$ is denoted by $\# V_B$.
 
   \begin{lem1} Given an orthonormal basis $\{e_i \}_{i=1}^n$ of $\G_n$, an $r$-vector $B\in \G_n^r$, expanded in this basis, is given by
   	\beq B =\sum_m B_{J_m}e_{J_m} ,    \label{plem}    \eeq
   	where $B_{J_m}=B\cdot e_{J_m}^\dagger$.
   \end{lem1}
     {\bf Proof:}    
   In the orthonormal $r$-vector basis $\{e_{J_m}\}$, the $r$-vector $B\in \G_n^r$
  is given by
  \[ B=\sum_{p=1}^{\pmatrix{n \cr r}} B_{J_p} e_{J_p} \]
  Dotting both sides of this equations with $e_{J_m}^\dagger$, the {\it reverse} of $e_{J_m}$, shows that
  $B_{J_m}=B\cdot e_{J_m}^\dagger$, since $e_{J_p}\cdot e_{J_m}^\dagger = \delta_{p,m}$.
  \smallskip
  
     $  \hfill \square$

\begin{def1}A non-zero $r$-vector $B\in \G_n^r$ is said to be divisible by a non-zero $k$-blade $K\in \G_n^k$, where $k\le r$, if $K B = \langle KB \rangle_{r-k}$. If
	$k=r$, $B$ is said to be totally decomposable. 
\end{def1}

 If $B$ is divisible by $L$, then 
 \[ B = \frac{L(LB)}{L^2}=L \frac{L\cdot B}{L^2}.     \]
 If $B$ is totally decomposable by $L$, then $B=\beta L$ for $\beta= \frac{L\cdot B}{L^2} \in \R$, in which case $B$ is an $r$-blade.

     A classical result is that an $r$-vector $B= \sum_{m} B_{J_m}e_{J_m}\in \G_n^r $ is an $r$-blade iff $B$ satisfies the {\it Pl\"ucker relations}. The Pl\"ucker relations, \cite[p.3]{poland}, are
   \beq (A\cdot B)\w B = 0  \label{pluckerrel} \eeq
   for all $(r-1)$-vectors $A \in \G_n^{r-1}$. In particular, for all coordinate $(r-1)$-vectors $e_{K_p}\in \G_n^{r-1}$, 
   \[  (e_{K_p}\cdot B)\w B=0  .    \]  
Since an $r$-blade trivially satisfies the Pl\"ucker relations, it follows that if there is a coordinate $(r-1)$-blade $e_{K}$, such that $(e_K \cdot B)\w B \ne 0$, then $B$ is not an $r$-blade.

For every coordinate $(r-1)$-blade $e_{K_p}$, there is a coordinate $r$-blade $e_{J_p}$, generally not unique, and a coordinate vector $e_{j_p}$, satisfying
\beq e_{J_p}=\pm e_{j_p}e_{K_p}\quad \iff \quad  e_{j_p} e_{J_p}=\pm e_{K_p} .   \label{ejcomment} \eeq
Indeed, when $e_{K_p}\cdot B \ne 0$, we always pick the coordinate vector
$e_{j_p}$, and rearrange the order of $K_p=\pm K_{j_p}$, in such a way that $B_{J_{j_p}}=(e_{j_p}\w e_{K_{j_p}})^\dagger \cdot B \ne 0$.
It follows that the Pl\"ucker relations (\ref{pluckerrel}) are equivalent to
   \[  e_{j_p}\cdot \Big[\Big( e_{K_{j_p}}\cdot  B\Big)\w B\Big]=B_{J_{j_p}}B- \Big[\Big(e_{K_{j_p}}\cdot B \Big)\w \Big(e_{j_p} \cdot B\Big) \Big]=0,           \]
   or equivalently,
   \beq  B =B_{J_{j_p}}^{-1}  \Big[\Big(e_{K_{j_p}}\cdot  B\Big)\w \Big(e_{j_p} \cdot B\Big) \Big],  \label{pluckerrel2} \eeq     
for all coordinate $(r-1)$-blades $e_{K_p} \in \G_n^{r-1}$ for which
$e_{K_p}\cdot B \ne 0$. It follows that when the Pl\"ucker relation (\ref{pluckerrel}) is satisfied, then the vector $e_{K_p}\cdot B$ divides $B$. 

    For a given non-zero $r$-vector $B=\sum_m B_{J_m}e_{J_m} \in \G_n^r$, let
    $\#_{max} S$ be the maximal number of linearly independent vectors in the set
    \beq S_B=\{ e_{K_j}\cdot B| \quad e_{K_j}\in \G_n^{r-1}  \}. \label{lem2} \eeq
    
    \begin{lem2} Given an $r$-vector $B\ne 0$, and the set $S_B$ defined in (\ref{lem2}).
    Then $r \le \#_{max}S \le n$, and $\#_{max}= r$ only if $B$ is an $r$-blade\end{lem2}  
  
    {\bf Proof:} Since $B\ne 0$, it follows that $B_{J_m}\ne 0$ for some 
  \[1\le m \le \pmatrix{n \cr r} .\]
  For each such $e_{J_m}$, there are $r$ possible choices for $e_{K_j}$ obtained by picking $K_j \subset J_m$, and the subset of $S_B$ generated by them are
  linearly independent. It follows that $r \le \#_{max}S_B \le n$. 
  
  Suppose now that $r < \#_{max}S_B$. To complete the proof of the Lemma, we show that there is a Pl\"ucker relation that is not satisfied for $B$, so it cannot be an $r$-blade. Let   \[   L=(e_{K_1}\cdot B)\w \cdots \w (e_{K_k}\cdot B)\ne 0,  \]
  be the outer product of the maximal number of linearly independent vectors from
  $S_B$. Then $k >r$ and $L B = L\cdot B \ne 0$. 
  
  Letting 
  \[ w = \big[(e_{K_1}\cdot B)\w \cdots \w (e_{K_{r+1}}\cdot B)\big]\cdot B \in \G_n^{1} , \]
  it follows that $w \ne 0$, and 
  \[ w^2 =  \big[(e_{K_1}\cdot B)\w \cdots \w (e_{K_{r+1}}\cdot B)\big]\cdot  (B \w w ) \ne 0.\]
   Since $w^2\ne 0$, it follows that $w\w  B \ne 0$. Noting the identity
   \[ w= \big[(e_{K_1}\cdot B)\w \cdots \w (e_{K_{r+1}}\cdot B)\big]\cdot B     \]
   \[ =(e_{K_1}\cdot B)  \big[(e_{K_2}\cdot B)\w \cdots \w (e_{K_{r+1}}\cdot B)\big]\cdot B \] 
   \[-(e_{K_2}\cdot B)  \big[(e_{K_1}\cdot B)\w \cdots  \lor^2 \cdots  \w (e_{K_{r+1}}\cdot B)\big]\cdot B  \]
   \[ + \cdots + (-1)^r (e_{K_{r+1}}\cdot B)  \big[(e_{K_1}\cdot B)\w \cdots \w (e_{K_r}\cdot B)\big]\cdot B ,   \]             
    it follows that
   \[ w \w B = \sum_{i=1}^r (-1)^{i+1}\alpha_i (e_{K_i}\cdot B)\w B \ne 0 , \]
   where $\alpha_i :=  \big[(e_{K_1}\cdot B)\w \cdots  \lor^i \cdots  \w (e_{K_{r+1}}\cdot B)\big]\cdot B\in \R  $.
  But this implies that
    $ \alpha_j (e_{K_j}\cdot B)\w B  \ne 0$ for some $j$, so the proof of the Lemma is complete.   
         
     $  \hfill \square$ 
     \smallskip
     
    It follows that if $\#_{max}S=r$, then 
   \[\alpha B= (e_{K_1}\cdot B)\w \cdots \w (e_{K_{r}}\cdot B). \]
   Furthermore, using (\ref{ejcomment}) and (\ref{pluckerrel2}), we can find distinct coordinate vectors $e_{j_1}, \cdots e_{j_r}$,
   and coordinate $(r-1)$-vectors $e_{K_{j_i}}=\pm e_{K_i}$, for $i=1, \cdots, r$, such
   that
   \[ \alpha = B e_{j_r \cdots j_1}=B\cdot  e_{j_r \cdots j_1} = B_{J_{j_2}}\cdots B_{J_{j_r}}.  \]   
     
    \begin{thm1} For $2 \le r \le n$, a non-zero $r$-vector $B\in \G_n^r$  is an $r$-blade iff the Pl\"ucker relations (\ref{pluckerrel}) are satisfied. 
    	When the Pl\"ucker relations are satified, then
    	\[  B =   B_{J_{j_2}}^{-1}\cdots B_{J_{j_r}}^{-1} (e_{K_{j_1}}\cdot B)\w \cdots \w (e_{K_{j _r}}\cdot B). \] 
    	    \end{thm1} 
      
      {\bf Proof:} 
      
      The proof follows directly from previous comments, and Lemma 2.
      
      \smallskip
      $  \hfill \square$
                   
    \section{Examples}  
      
    Everybodies' favorite example is when $r=2$ and $n=4$. A general non-zero $2$-vector $B\in \G_4^2$ has the form  
    \[ B = \sum_{m=1}^6 B_{J_m}e_{J_m}\in \G_4^2. \]
      Since $B\ne 0$, at least one component $B_{J_p}e_{ij}\ne 0$. Without loss of
    generality, we can assume $e_{ij}=e_{12}$. It follows that the Pl\"ucker relation (\ref{pluckerrel2}) for this component is
    \[  B=B_{12}^{-1}\Big[ \big(e_1\cdot B\big)\w \big(e_2\cdot B\big)  \Big]  \]
    \beq  =B_{12}^{-1}\Big[ (B_{12}e_2 +B_{13} e_3 +B_{14}e_4 )\w (-B_{12}e_1+B_{23}e_3+B_{24}e_{4}   )  \Big].   \label{favorex} \eeq
      
    Since
    \[ 2 (e_i \cdot B)\w B= e_i\cdot (B\w B)  =(B_{12}B_{34}-B_{13}B_{24}+B_{14}B_{23}) e_ie_{1234},      \]
   it follows that $B$ is a $2$-blade iff $B\w B=0$, or equivalently,
    \beq B_{12}B_{34}-B_{13}B_{24}+B_{14}B_{23} =0. \label{pluck3}  \eeq
     This occurs when $B$ is totally decomposable, given in (\ref{favorex}).
    More generally, a similar argument holds true for $r=2$ and any $n>4$, but with more Pl\"ucker relations (\ref{pluck3}) to be satisfied.
    
   To get more insight into what is going on, let $B=\sum_{m} B_{J_m}e_{J_m}\in \G_n^r$ be a general $r$-vector.
 Following \cite{d-b-N}, we calculate
 \[ B^2 = \sum_{m,p} B_{J_m}B_{J_p}e_{J_m}e_{J_p} = \sum_{m = p}   B_{J_m}B_{J_p}e_{J_m}e_{J_p}+ \sum_{m\ne p} B_{J_m}B_{J_p}e_{J_m}e_{J_p}      \] 
 \beq  =  \sum_{m = p}   B_{J_m}B_{J_p}e_{J_m}e_{J_p}+ \sum_{m< p} B_{J_m}B_{J_p}\Big(1 + (-1)^{r-k}\Big)e_{J_m}e_{J_p} ,   \label{calBB} \eeq
 where $k=\# J_m \cap J_p$, showing that $B^2 = \langle B^2 \rangle_0$ if $r$ and $k$ have different parity for all values of $m<p$.  
 
 Consider $B=e_{123}+e_{456} \in \G_6^3$. It is easy to show, using (\ref{pluckerrel}),
 that $B$ is not a $3$-blade, since
 \beq (e_{12}\cdot B)\w B = e_3 \w B =e_{3456} \ne 0.    \label{no3blade} \eeq
Note, however, that $B^2 = -2$, as also follows from (\ref{calBB}) since $k$ and $r$ have opposite parities. Clearly, the condition $B^2-\langle B^2 \rangle_0 =0$ is not sufficient to guarantee that $B$ is an $r$-blade, but $B^2 \ne \langle B^2 \rangle_0 $ 
guarantees that $B$ is {\bf not} an $r$-blade.
   
 In \cite{d-b-N}, Nguyen expresses the Pl\"ucker Relations differently, and gives a  different proof. The following Theorem shows that Nguyen's Pl\"ucker Relations are
 equivalent to our definition (\ref{pluckerrel}). Let $e_K$ denote an arbitrary coordinate $(r-1)$-vector in $\G_n^{r-1}$.

\begin{thm2} A non-zero $r$-vector $B\in \G_n^r$ is an $r$-blade iff 
	both the conditions
	\[  B^2 = B\cdot B, \ \ {\rm and} \ \  BvB \in \G_n^1 \ \ {\rm for}  \ {\rm all} \ \ v=e_K \cdot B \in \G_n^1,      \]
	are satisfied.	
\end{thm2}
 
   {\bf Proof:} It is easy to show that if $B \in \G_n^r$ is an $r$-blade, then both of the 
 conditions in the Theorem are true. To complete the proof we must only show that if both conditions are not satisfied, then $B$ is not an $r$-blade. If $B^2 \ne B\cdot B$, then $B$ is not an $r$-blade, so it is only left to show that if $B^2=B\cdot B$ and $BvB\not \in \G_n^n$, then $B$ is not an $r$-blade. 
  
  The following identity is needed:
  \[ B v B = B (v\cdot B + v\w B)= B(v\cdot B)+B(v\w B)  \] 
 \beq =B\cdot (v\cdot B) + B \cdot(v\w B)+  \langle  B (v B)\rangle_{3,5,\cdots}  \label{ident} \eeq
 Suppose now that $B^2 = B\cdot B$ and $BvB \not \in \G_n^1$ for some $v=e_K \cdot B$. The identity (\ref{ident}) implies that $\langle BvB \rangle_{3,5,\cdots}\ne 0$. We also know that
 \[ 2B\w v = Bv + (-1)^r vB\quad \iff \quad Bv = 2B\w v -(-1)^r vB,  \]
 so that
 \[ \langle BvB \rangle_{3,5,\cdots}= \langle \big( 2B\w v -(-1)^r vB  \big)B \rangle_{3,5,\cdots} =  2\langle \big( B\w v  \big)B \rangle_{3,5,\cdots}\ne 0,  \]
 since $B^2=B\cdot B$. But this implies that $ B\w (e_K\cdot B) \ne 0$, so by
 (\ref{pluckerrel}), $B$ is not an $r$-blade.
 $  \hfill \square$

\smallskip

\begin{figure}[h]
	\begin{center}
		\includegraphics[width=0.2 \linewidth]{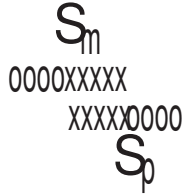}
		\caption{Shows $S_m$ and $S_p$ in the case for odd $r=9$, with an odd overlap $k=5$. In this case we choose $e_K=e_{j_2 \cdots j_9}$, the $j_2,\cdots j_9$ are the last $8$ digits of $S_m$. A similar construction applies when $r$ and $k$ are even.}
		\label{pluckerfig}
	\end{center}
\end{figure}

Note that (\ref{pluckerrel}) is fully equivalent to the conditions in the Theorem. 
If $B^2 \ne B\cdot B$ then there is at least one $e_K \in \G_1^{r-1}$ such that $(e_K \cdot B)\w B \ne 0$. This follows from the coordinate expansion of $B$ given in (\ref{calBB}), because the parity condition allows us to find at least one non-zero $term(m,p):=B_{J_m}B_{J_p}e_{J_m}e_{J_p}$, with overlap $k=\# (S_{J_m}\cap S_{J_p})$, such that
$(e_K\cdot B)\w B \ne 0$. We simply pick $e_K$ in such a way that $K \subset S_{J_m}$,
but $\# (K\cap J_p)\le r-2$, see Figure \ref{pluckerfig}. 

It is also interesting to note that just as $B^2=B\cdot B$ whenever $B$ is an
$r$-vector satisfying the Pl\"ucker Relations (\ref{pluckerrel}),
$BvB=(-1)^{r+1} (B\cdot B)v$ for all $v=e_K\cdot B$ whenever $B$ is an $r$-vector satisfying the Pl\"ucker Relations (\ref{pluckerrel}). This follows from the identity
\[ BvB = B[2 v\w B+ (-1)^{r+1} Bv]=(-1)^{r+1} B^2v \quad \iff \quad v \w B = 0 \]
for all $v=e_K\cdot B$, whenever $B$ satisfies the Pl\"ucker Relations (\ref{pluckerrel}).

An example of a decomposable $3$-vector $B\in \G_{5}^3$  is
\[ B= e_{125} +e_{234}+2 e_{124}+e_{235}+e_{123}+e_{245},  \]
satisfying the condition
\[  B_{125}B_{234}-B_{124}B_{235}+B_{123}B_{245}=0.  \]
One factorization is $B=$
\[ (-e_1+e_2+e_3+e_4)\w (-2e_1+e_2+e_3-e_5)\w (-3e_1+e_2+2e_3+e_4-e_5).     \]

To gain further insight, consider $B\in \G_3^6$ and where $B^2 =\langle B^2 \rangle_0$.
For example, let
\[  B=e_{123}+e_{456} +  e_{124}+e_{356} +e_{125}+e_{346}+  e_{126}+e_{345}.   \]
Then,
\[ (e_{21}\cdot B)\w B =(B_{123}B_{456} -  B_{124}B_{356} +B_{125}B_{346}-  B_{126}B_{345})e_{3456}     \]
\[+( B_{123}B_{145}-B_{124}B_{135}+B_{125}B_{134})e_{1345} + \cdots  \] 
\[ +( B_{124}B_{156}-B_{125}B_{146}+B_{126}B_{145})e_{1456}=0,   \]
but
\[ (e_{13}\cdot B)\w B = e_2\w B = e_{2456}+e_{2356}+e_{2346}+e_{2345}\ne 0.   \]
Since $B$ satisfies the first condition, it is divisible by $e_{12}\cdot B$. But since it does not satisfy the second condition, it is not totally decomposable.

\section*{Acknowledgement} 
I thank Dr. Dung B. Nguyen for calling my attention to the problem of expressing and proving the Pl\"ucker relations in geometric algebra. This work is the result of many email exchanges with him, false starts, and trying to understand the roots of the vast literature on the subject and its generalizations.


\begin{thebibliography}{}
 \bibitem{ardila}F. Ardila, {\it The Geometry of Matroids}, Notices of the AMS, Vol. 65, No. 8, Sept. 2018.
 \bibitem{dandw91}A. Dress and W. Wenzel, {\it Grassmann-Pl\"ucker Relations and Matroids with Coefficients}, Advances in Mathematics 86, 68-110 (1991).
 	\bibitem{FH1991}W. Fulton, J. Harris, {\it Representation Theory:
 		A First Course}, Springer-Verlag 1991.
 	\bibitem{harris1992}J. Harris, {\it Algebraic Geometry: A First Course}, Spring, p.63-71, 1992.
 	\begin{verbatim}http://math.rice.edu/~evanmb/math465spring11/\end{verbatim}
 	\begin{verbatim}Harris-Grassmannians.pdf\end{verbatim} 
 	\bibitem{H/S} D. Hestenes and G. Sobczyk. {\it Clifford Algebra to
 			 		Geometric Calculus: A Unified Language for Mathematics and Physics},
 				2nd edition, Kluwer 1992.
 	
 	\bibitem{milne}J. Milne, {\it Algebraic Geometry}, Version 6.06, March 2017, pp. 144-148. 
 	http://www.jmilne.org/math/CourseNotes/AG.pdf
 	
 	\bibitem{d-b-N}D. B. Nguyen, {\it Pl\"ucker's Relations, Geometric Algebra, and the Electromagnetic Field}, Version: Sept. 3, 2018.
 	
 	\bibitem{geoS2017} G. Sobczyk, {\it Geometrization of the Real Number System}, July, 2017.
 	\begin{verbatim}http://www.garretstar.com/geonum2017.pdf\end{verbatim}
 	\bibitem{poland} \begin{verbatim}
 	http://homepages.math.uic.edu/~coskun/poland-lec1.pdf
 	\end{verbatim}
 	\bibitem{SNF} G. Sobczyk, {\em New Foundations in Mathematics: The Geometric Concept of Number},
 	\newblock Birkh\"auser, New York 2013.		
  \bibitem{ctc-wall2008}C.T.C. Wall, {\it Pl\"ucker formulae for curves in high dimensions},
University of Liverpool, 2008.
\begin{verbatim}
https://www.liverpool.ac.uk/~ctcw/pluckD.pdf
\end{verbatim}
 
 \end{thebibliography}
\end{document}